\documentclass[12pt]{article} \usepackage{amssymb,amsfonts,epsf}
\catcode`\@=11\renewcommand{\@date}{30 June 2000}
\newcommand{\vs}{\vspace{10pt}}\newcommand{\y}{\\[3pt]}
\newcommand{\q}{\quad}\newcommand{\m}{P}\renewcommand{\=}{\!=\!}
\renewcommand\ge{\geqslant}\renewcommand\le{\leqslant}
\newcommand\bb{\bigbreak}
\newcommand\mb{\medbreak}\renewcommand\sb{\smallbreak}
\newcommand{\bs}{\backslash}\newcommand{\g}{\mathfrak{g}}
\newcommand\n{\noindent}\newcommand{\bu}{\hbox{\large$\bullet$}\q}

\newcommand{\cC}{\mathcal{C}}
\newcommand{\ZZ}[1]{\al#1\ar^\perp\kern-8pt\lower2pt\hbox{$_\cZ$}}
\newcommand{\EE}[1]{\lceil#1\kern-0.5pt\rfloor}
\newcommand{\cE}{\mathcal{E}}\newcommand{\cF}{\mathcal{F}}
\newcommand{\cS}{\mathcal{S}}\newcommand{\cZ}{\mathcal{Z}}
\newcommand{\D}{\mathbb{D}}\newcommand{\id}{I}
\newcommand{\C}{\mathbb{C}}\newcommand{\R}{\mathbb{R}}
\newcommand{\CP}{\mathbb{CP}}\newcommand{\W}{\mathcal{W}}
\newcommand{\GH}{G_\mathrm{H}}\newcommand{\gH}{\g_\mathrm{H}}
\newcommand{\suml}[2]{\hbox{$\textstyle\sum\limits_{#1}^{#2}$}}
\newcommand{\dm}{d\kern1pt\raise.8pt\hbox{$^\m$}}
\newcommand\CS{C}
\newcommand\si{\sigma}\newcommand\la{\lambda}
\renewcommand{\th}{\theta}\newcommand\ep{c}
\newcommand\w{\omega}\newcommand\Ga{\Gamma}
\newcommand{\ph}{\phantom}\newcommand{\we}{\wedge}

\renewcommand{\ll}{[\![}\newcommand{\rr}{]\!]}
\newcommand{\op}{\oplus}\newcommand{\ot}{\otimes}
\newcommand\qed{\hfill$\square$\mb}\newcommand{\rf}[1]{(\ref{#1})}
\newcommand{\E}{\raise1pt\hbox{$\textstyle\bigwedge$}\kern-0.5pt}
\newcommand{\ba}{\begin{array}}\newcommand{\ea}{\end{array}}
\newcommand{\be}{\begin{equation}}
\newcommand{\ee}[1]{\label{#1}\end{equation}}
\newcommand{\bt}{\begin{tabular}}\newcommand{\et}{\end{tabular}}
\newcommand{\al}{\langle}\newcommand{\ar}{\rangle}
\newcommand{\nbf}[1]{\bigbreak\n\textbf{#1.}}
\newcommand{\nit}[1]{\sb\n\textit{#1.}}

\newcommand{\frt}[2]{\hbox{$\textstyle\frac{#1}{#2}$}}

\newcommand{\z}{\\[-3pt]\small}
\begin{document}

\title{\bf Almost Hermitian Geometry\y on Six Dimensional Nilmanifolds}\vs\vs

\author{E.~Abbena\thanks{Work supported by the CNR research group GNSAGA and
the funds of MURST.}\z Dipartimento di Matematica\z Universit\`a di Torino \and
S.~Garbiero$^*$\z Dipartimento di Matematica\z Universit\`a di Torino \and
S.~Salamon\z Mathematical Institute\z University of Oxford}\maketitle
 
\nit{Abstract} The fundamental 2--form of an invariant almost Hermitian
structure on a 6--dimensional Lie group is described in terms of an action by
$SO(4)\times U(1)$ on complex projective 3--space. This leads to a
combinatorial description of the classes of almost Hermitian structures on the
Iwasawa and other nilmanifolds.\mb

\nit{AMS classification} 53C15; 53C55, 51A05, 17B30.

\subsection*{Introduction}

Let $(M,g,J)$ be an almost Hermitian manifold, so that $J$ is an almost
complex structure orthogonal relative to the Riemannian metric $g$. The
decomposition of $\nabla J$ ($\nabla$ is the Levi--Civita connection) into
irreducible components under the action of the unitary group determines the
Gray--Hervella class of the almost Hermitian structure \cite{GH}. In real
dimension $2n$ with $n\ge3$, $\nabla J$ has four components $W_1,W_2,W_3,W_4$,
and one is interested in structures for which one or more of these
vanishes. For example, $(M,J)$ is a complex manifold if and only if
$W_1=W_2=0$, and $(M,\w)$ is a symplectic manifold ($\w$ is the fundamental
2--form determined by $g$ and $J$) if and only if $W_1=W_3=W_4=0$.

In this paper, we are concerned with left--invariant tensors $g,J$ on a Lie
group $G$ of real dimension 6. Indeed, we fix a metric $g$ and consider the
space $\cZ$ of \textit{all} left--invariant almost complex structures $J$
compatible with $g$ and an orientation. This reduces many questions to
properties of the Lie algebra $\g$ of $G$. The manifold $\cZ$ is isomorphic to
the complex projective space $\CP^3$, and a choice of standard coordinates
allows us to visualize it in terms of a tetrahedron, in which the edges and
faces represent projective subspaces $\CP^1$ and $\CP^2$. This technique and
was first used in \cite{AGS}, which identified the subsets $\cC$, $\cS$ of
$\cZ$ corresponding to complex and symplectic structures compatible with a
given metric on the complex Heisenberg group $\GH$.

In order to develop the theory in a more systematic way, we show in \S3 that
every single component of $\nabla J$ (equivalently, $\nabla\w$) can be readily
extracted by means of wedging with appropriate differential forms. Some of the
integrability equations (such as those of Lemma~2) interact effectively with
the nilpotency condition on \be d:\g^*\to\E^2\g^*\ee{d} that arises in the
theory of minimal models, and our approach is particularly suited to the case
in which the Lie group $G$ is nilpotent. On the other hand, our investigation
is also designed to illustrate aspects of the theory of differential forms and
complex structures on 6--manifolds that can be applied to non--invariant
settings.

The non--existence of a K\"ahler metric on a nilmanifold $\Ga\bs G$ (other than
a torus) will imply that there are no compatible complex and symplectic
structures, so $\mathcal{C}\cap\mathcal{S}=\emptyset$. This relatively deep
fact pervades the descriptions of subsets of $\cZ$ in our examples, and
renders the combinatorial aspect of the results all the more striking. In six
dimensions, a fundamental 2--form $\w$ satisfies $W_4=0$ if and only if
$\w\we\w$ is closed, and this is an especially natural condition in our
context. We point out that the set of such 2--forms is orthogonal to the image
of \rf{d} and always non--empty. The corresponding class of cosymplectic
structures is an intersection of real hypersurfaces of $\cZ$, and typically
has dimension $b_1$ (the first Betti number of $\g$).

The kernel of \rf{d} gives rise to a real 4--dimensional space $\D$ of
invariant closed 1--forms on $\GH$, and the action of a corresponding subgroup
$SO(4)$ on $\cZ$ provides symmetry that can be exploited to simplify the
equations. In fact, we parametrize $\w(\m;a,b)\in\cZ$ by means of $\m\in
SO(4)$ and $(a,b)$ in a unit circle, and seek solutions $\m;a,b$ to the
problem at hand, exploiting the concept of self--duality and a `conjugated'
exterior derivative $\dm$. Most of the calulations were done by hand, and then
checked with \textsc{Maple}'s differential form package. The most interesting
part is to interpret the conditions geometrically as subsets of $\cZ$, and
there is a sufficiently rich class of examples to illustrate many contrasting
features of the theory.

Once one has determined the four classes for which $W_i$ vanishes (with
$i=1,2,3,4$ in turn), all the other classes may be obtained by taking
intersections. The resulting inclusions can be summarized by means of a
quotient lattice, which is illustrated in \S4 for $G_H$. In this case, we
obtain a complete description of the 16 Gray--Hervella classes, in terms of
faces, edges, and vertices of the standard tetrahedron. The upshot is that
every class (for which at least one $W_i$ vanishes) is contained in the union
of two particular faces, and has at most two connected components.

The techniques are equally applicable to a class of similar examples, and we
carry out the same process for the two other irreducible nilpotent Lie
algebras with $b_1=4$ (indeed, $b_1\ge4$). In these cases, we give a complete
description of the classes for which one of the `larger' components $W_2,W_3$
is zero, and these classes get progressively smaller in the three cases. We
also identify enough subsets (unions of points, circles and 2-spheres) with
$W_1=0$ or $W_4=0$ to determine most of the classes for which two $W_i$
vanish. In \S5, we show that the other 2--step example has exactly four
Hermitian structures, but that two of these can be made to coincide by
modifying the inner product. In \S6, we are able to detect the non--existence
of both complex and symplectic structures in the 3--step case, for which the
various classes are best visualized using a `similar' tetrahedron arising from
a change of basis.

\subsection*{1.~Invariant tensors in 6 dimensions}

The complex Heisenberg group is given by the set of matrices \be
\GH=\left\{\left(\ba{ccc}1&\zeta_1&\zeta_3\\0&1&\zeta_2\\0&0&1\ea\right):
\>\zeta_i\in\C,\ i=1,2,3\right\}\ee{GH} under multiplication. The $1$--forms
$\alpha^1=d\zeta_1$, $\alpha^2=d\zeta_2$, $\alpha^3=-d\zeta_3+\zeta_1d\zeta_2$
are left--invariant on $\GH$. The Lie bracket is determined by the equations
\be\left\{\ba{l} d\alpha^i=0,\q i=1,2,\y d\alpha^3=\alpha^1 \we\alpha^2,
\ea\right.\ee{dalpha} and the $\alpha^i$ are holomorphic relative to the
natural complex structure $J_0$ on $\GH$. Mapping the above matrix to
$(\zeta_1,\zeta_2)$ determines a homomorphism $\mu\colon G_H\to(\C^2,+)$ of
complex Lie groups.

In this paper, we shall treat $\GH$ as a \textit{real} Lie group, and setting
\[\alpha^1 =e^1+ie^2,\q \alpha^2=e^3+ie^4,\q \alpha^3=e^5+ie^6\]
provides a real basis $(e^i)$ of $\gH^*$ where $\gH$ denotes the corresponding
Lie algebra. The real version of \rf{dalpha} is therefore \be\left\{\ba{l}
de^i=0,\q 1\le i\le4,\y de^5=e^1\we e^3+e^4\we e^2\y de^6=e^1\we e^4+e^2\we
e^3.\ea\right.\ee{rv} The kernel of \rf{d} is the 4--dimensional subspace $\D$
spanned by $e^1,e^2,e^3,e^4$, and this coincides with the image of $\mu^*$.
Moreover, the image of \rf{d} lies in $\E^2\D$. The fundamental role played by
$\D$ in the theory of invariant structures on $\GH$ was emphasized in
\cite{AGS}, and we now set about generalizing this situation.

Let $G$ be a real $6$--dimensional nilpotent Lie group, and let $\Ga$ be a
discrete subgroup of $G$ for which the set $\Ga\bs G$ of right cosets is a
compact manifold $M$. A classification of corresponding Lie algebras was given
in \cite{Mag} (see the tables in \cite{CFGU,S}). Given that $G$ necessarily
has rational structure constants, such a $\Ga$ must exist \cite{Mal}.  For
example, the \textit{Iwasawa manifold} is the compact quotient space $M=\Ga\bs
\GH$ where $\Ga$ is the subgroup defined by restricting $\zeta_i$ in \rf{GH}
to be Gaussian integers. The existence of a compact quotient enables one to
apply Nomizu's theorem to compute the cohomology of $\Ga\bs G$, and to assert
the non--existence of a K\"ahler metric on $M$ unless $G$ is abelian (see
\cite{H,BG,McD,C,TO}).

We shall suppose that a 4--dimensional subspace $\D$ of $\g^*$, the dual of the
Lie algebra of $G$, is chosen. Let $(e^1,\dots,e^6)$ a basis of $\g^*$ such
that $\D=\al e^1,e^2,e^3,e^4\ar$. Orientations on $\g^*$ and $\D$ are
determined by the forms \be\upsilon=e^{123456},\q\upsilon'=e^{1234} \ee{up}
respectively. (From now on we abbreviate $e^j\we e^j\we\cdots$ to
$e^{ij\cdots}$.) Each $e^i$ is a left--invariant 1--form on $G$, and therefore
passes to the quotient $M$. We shall denote the corresponding form on $M$ by
the same symbol; for example, $\upsilon$ determines an orientation of $M$. A
tensor on $M$ will be called \textit{invariant} if it can be expressed in
terms of the basis $(e^i)$ using constant coefficients. In particular, a
differential form is invariant if and only if its pull--back to $G$ is
left--invariant by $G$. Declaring the chosen $1$--forms $e^i$ to be
orthonormal determines an invariant Riemannian metric \be g=\suml{i=1}6
e^i\otimes e^i.\ee{g} In the following study, we shall regard the choice
of $\upsilon$, $\D$, $\upsilon'$, $g$ as fixed, rather than a particular basis
$(e^i)$.

The symmetry group of $(\g,\upsilon,g)$ is $SO(6)$ and that of $(\D,\upsilon',
g)$ is $SO(4)$. The decomposition \be\E^2\D=\E^2_+\D\op\E^2_-\D \ee{+-} into
eigenvalues of the $*$ operator gives rise to a double covering $SO(4)\to
SO(3)\times SO(3)$. It follows that $\m\in SO(4)$ can be represented by the
$6\times6$ matrix \be\left(\ba{cc}\m_+&0\\0&\m_-\ea\right),\ee{SO4} where
$\m_+:\E^2_+\D\to\E^2_+\D$ and $\m_-:\E^2_-\D\to\E^2_-\D$. With respect to the
bases \[\ba{rcl}(e^{12}+e^{34},e^{13}+e^{42},e^{14}+e^{23})&\hbox{of}&\E^2_+
\D,\y (e^{12}-e^{34},e^{13}-e^{42},e^{14}-e^{23})&\hbox{of}&\E^2_-\D,\ea\] we
may write \be \m_+=\left(\ba{ccc}r&u&x\\s&v&y\\t&w&z\ea\right),\q
\m_-=\left(\ba{ccc}r'&u'&x'\\s'&v'&y'\\t'&w'&z'\ea\right)\in SO(3).\ee{m+-}

An invariant almost complex structure $J$ on $M$ is an endomorphism of $\g$
(or of $\g^*$) such that $J^2=-1$. It is uniquely determined by writing
\[Je^j=\suml{i=1}6 a_i^je^i,\q j=1,\dots,6\] where $(a_i^j)$ is a constant
matrix. Following the notation of \cite{AGS}, we denote by $\cZ$ the set of
invariant positively--oriented orthogonal almost complex structures on
$M$. Thus, $J\in\cZ$ if and only if $(a_i^j)$ is both orthogonal and
skew--symmetric with positive determinant. The group $SO(6)$ acts on $\cZ$ by
conjugation, and the stabilizer of the almost complex structure $J_0$, for
which (like the one defined by \rf{dalpha}) \be J_0e^1=-e^2,\q J_0e^3=-e^4,\q
J_0e^5=-e^6,\ee{J0} can be identified with $U(3)$. Thus, $\cZ$ is the
Hermitian symmetric space $SO(6)/U(3)$. The double covering $SU(4)\to SO(6)$
allows one to regard $SO(6)$ as a transitive subgroup of projective
transformations of $\CP^3$, and the following is well known.

\nbf{Proposition 1} $\cZ$ is isomorphic to the complex projective space
$\CP^3$.\mb 

We describe the isomorphism explicitly, referring the reader to \cite{AGS,Ap}
for more details. Let $V\cong\C^4$ denote the standard representation of
$SU(4)$ and let $(v^0,v^1,v^2,v^3)$ be a unitary basis of $V$. Then $\E^2V$ is
the complexification of a real vector space which we identify with
$\g^*$. This identification can be chosen in such a way that \be\ba{ll}
2v^{01}=e^1+ie^2, \qquad & 2v^{23}=e^1-ie^2,\y 2v^{02}=e^3+ie^4,\qquad &
2v^{31}=e^3-ie^4,\y 2v^{03}=e^5+ie^6,\qquad & 2v^{12}=e^5-ie^6\ea\ee{vij}
(remember that $v^{ij}$ denotes $v^i\we v^j$). A point $J\in\cZ$ corresponds
to a totally isotropic subspace of the complexification of $\g^*$, namely the
$i$--eigenspace of $J$, and any such subspace equals \[V_u=\{u\we v:v\in V\}
\subset\E^2V,\] which depends on $[u]\in\mathbb{P}(V)\cong\CP^3$. For example,
\rf{J0} corresponds to the point $[v^0]=[1,0,0,0]$ of $\CP^3$. We shall also
consider the almost complex structures \[J_1=[0,1,0,0],\q J_2=[0,0,1,0],\q
J_3=[0,0,0,1].\]

Notice that an almost complex structure with the `wrong' orientation
corresponds not to a subspace of type $V_u$, but to one of type $\E^2W$ where
$W$ is a 3--dimensional subspace of $V$. For example, $-J_0$ corresponds to
$W=\al v^1,v^2,v^3\ar$. This is an aspect of the well-known
$\alpha\leftrightarrow\beta$ duality of the Klein correspondence
\cite{Penrose}. In fact many of the constructions below illustrate concepts
from elementary projective geometry.

\subsection*{2.~An action of $SO(4)$ on $\CP^3$}

The \textit{fundamental 2--form} $\w$ of the almost Hermitian structure
$(g,J)$ is defined by \[\w(X,Y)=g(JX,Y).\] Associated to $J_0,J_1,J_2,J_3$ are
the fundamental 2--forms \be\ba{rcl} \bu\w_0
&=&\;e^{12}+e^{34}+e^{56},\\\bu\w_1&=&\;e^{12}-e^{34}-e^{56},\\\bu\w_2
&=&-e^{12}+e^{34}-e^{56},\\\bu\w_3 &=&-e^{12}-e^{34}+e^{56}.\ea\ee{vert} Since
the metric $g$ will always be fixed, one may equally well label points of
$\cZ$ by their corresponding fundamental forms.

Given the coordinate system on $\CP^3$ determined by $(v^i)$, it is helpful to
visualize $\cZ$ as a \textit{solid} tetrahedron with vertices \rf{vert}, and
to consider \\\ph{mmmmmmmm}\bt{l} $\cE_{ij}$, the `edge' ($\cong\CP^1\cong
S^2$) containing $\w_i,\w_j$;\y $\cF_i$, the `face' ($\cong\CP^2$) opposite
$\w_i$.\et\par\n Whilst this picture may be geometrically inaccurate, it
provides a powerful tool for analysing some of the singular varieties that
arise in the classification of almost Hermitian classes.

For future reference, we define a number of additional objects:\sb

\n(i) The \textit{equatorial circle} in the edge $\cE_{ij}$ is the set
\be\CS_{ij}=\{\w\in \cE_{ij}:g(\w,\w_i)=g(\w,\w_j)\}\ee{equat} where $g$ here
denotes the metric induced on 2--forms. For example, \be \ep(\th)=e^{56}
+\cos\th(e^{13}+e^{42})+\sin\th(e^{14}+e^{23})\ee{eth} is a typical element of
$\CS_{03}$.\sb

\n(ii) The \textit{generalized edge} determined by a decomposable unit 2--form
$\si=e\we f$ is \be\EE{\si}=\{\si+\tau\in\cZ: \tau\in\E^2\al
e,f\ar^\perp\}.\ee{EE} Consistency with the orientation (see \rf{up}) requires
that $\tau$ be chosen from a 2-sphere, and $\EE{\si}$ is a typical complex
projective line in $\cZ$. From \rf{vij}, $\cE_{03}=\EE{e^{56}}$ is the line
$\mathbb{P}(\al v^0,v^3\ar)$.\sb

\n(iii) The \textit{polar set} of an arbitrary non--zero 2--form $\si$ is
\[\ZZ{\si}=\{\w\in\cZ:g(\w,\si)=0\};\] this is simply the intersection of the
hyperplane $\al\si\ar^\perp$ of $\E^2\R^6$ with the submanifold $\cZ$ of
fundamental forms with respect to $g$.\sb

We shall now describe the action of $SO(4)=\hbox{Aut}(\D,\upsilon',g)$ on
$\cZ$. Let $J$ be an almost complex structure with fundamental 2--form
$\w\in\cZ$. Since $-Je^5$ is a unit 1--form orthogonal to $e^5$, it has the
form $ae^6+bf^1$ where $f^1\in\D$, $\|f^1\|=1$, and $a^2+b^2=1$. The unit
1--form $af^1-be^6$ is then orthogonal to both $e^5$ and $Je^5$. The next
result is a consequence of this observation.

\nbf{Proposition 2} The fundamental 2--form defined by an arbitrary point of
$\cZ$ has the form \be\w=e^5\we(ae^6+bf^1)-f^2\we(af^1-be^6)+f^3\we f^4,\ee{w}
where $(f^1,f^2,f^3,f^4)$ is an oriented orthonormal basis of $\D$ and
$a^2+b^2=1$.\bb

Given $(f^i)$, we may define $\m\in SO(4)$ by setting $\m(e^i)=f^i$ for $1\le
i\le4$. Relative to the original basis $(e^i)$, it therefore makes sense to
write \rf{w} as $\w(\m;a,b)$. For example, $\w(\id;1,0)=\w_0$ and
$\w(\id;-1,0)=\w_2$ ($\id$ denotes the identity) are two of the vertices of
the tetrahedron, and \be\w(\id;a,b)=e^{34}+\tau,\q\tau\in\E^2\al e^1,e^2,e^5,
e^6\ar\ee{greatc} lies in $\cE_{02}=\EE{e^{34}}$ for any $(a,b)\in S^1$. In
particular, the points \be\ba{l} \varpi_0=\w(\id,0,1)=-e^{15}+e^{26}+e^{34},\y
\varpi_2=\w(\id,0,-1)=e^{15}-e^{26}+e^{34}\ea\ee{varpi} lie in the circle
$\CS_{02}$.

What the Proposition is really saying is that there is a transitive action of
$K=SO(4)\times U(1)$ on $\cZ$, where $SO(4)$ is described in \rf{SO4}, and
$U(1)$ acts by rotating points $(a,b)$ of the circle. Replacing $b$ by $-b$ in
\rf{w} has the same effect as changing the sign of every $f^i$, so the element
$(-\id,-1)\in K$ acts trivially, and there is an effective action of
$K/\mathbb{Z}_2$ on $\cZ$. Consider the possible stabilizers of $SO(4)$ on
$\cZ$ and the corresponding orbits:

\nit{Exceptional case} If $b=0$ then $a=\pm1$, $\D$ is $J$--invariant and the
subgroup of $SO(4)$ stabilizing $\w(\m;\pm1,0)$ is isomorphic to $U(2)$. This
gives rise to two exceptional orbits: as $\m$ varies in $SO(4)$, \par\sb\bt{l}
$\w(\m;1,0)$ spans the edge $\cE_{03}=\EE{e^{56}}$, and\y $\w(\m;-1,0)$ spans
the edge $\cE_{12}=\EE{-e^{56}}$.\et\sb

\nit{Generic case} If $b\ne0$, the stabilizer of $\w(\m;a,b)$ is the subgroup
$SO(2)$ of $SO(4)$ fixing $f^1$, $f^2$, $f^{34}$, and acting by the rotation
\be\left\{\ba{rcl}f^3 &\mapsto& \cos\theta\,f^3+\sin\theta\,f^4,\y f^4
&\mapsto&\!\!-\sin\theta\,f^3+\cos\theta\,f^4.\ea\right.\ee{rot} The generic
orbit of $SO(4)$ on $\cZ$ is therefore the 5--dimensional space $SO(4)/SO(2)$.
An example of such an orbit is \be\ZZ{e^{56}}=\{\m(\varpi_0)=\w(\m;0,1):\m\in
SO(4)\};\ee{perp} this is because the points in \rf{varpi} are orthogonal to
the 2--form $e^{56}$, which is fixed by $SO(4)$.

\mb One can now visualize the action of $K$ filling out the whole tetrahedron,
and the relevance of \rf{m+-}. Each edge represents a complex projective line,
that is, a 2-sphere. As $a,b$ vary, \rf{w} describes a great circle on another
such line or 2-sphere $\Sigma$, corresponding to the edge shown bold in the
figure. We view it as an elongated elliptical wire, one point $p_+=\w(\m;1,0)$
of which lies on $\cE_{03}$ and the diametrically opposite point
$p_-=\w(\m;-1,0)$ on $\cE_{12}$. Varying $\m_+$ moves $p_+$ `up and down'
$\cE_{03}$ (`around' would be more accurate), and varying $\m_-$ moves $p_-$
around $\cE_{12}$.

The stabilizer of the pair $(p_+,p_-)$ in $SO(4)$ is isomorphic to $U(1)$, and
this acts to rotate the wire to fill out the 2-sphere $\Sigma$ with poles
$p_+,p_-$. This is what is represented by the right-hand side of the
Figure. The forms \rf{greatc} may be regarded as the wire's `position of
rest', in which the bold edge coincides with $\cE_{02}$. The label `$\pm ij$'
on an edge indicates that it equals $\EE{\pm e^{ij}}$ in the notation of
\rf{EE}.

\phantom.\vspace{160pt}
\phantom.\hspace{10pt}\epsfxsize150pt\epsfbox[0 0 205 160]{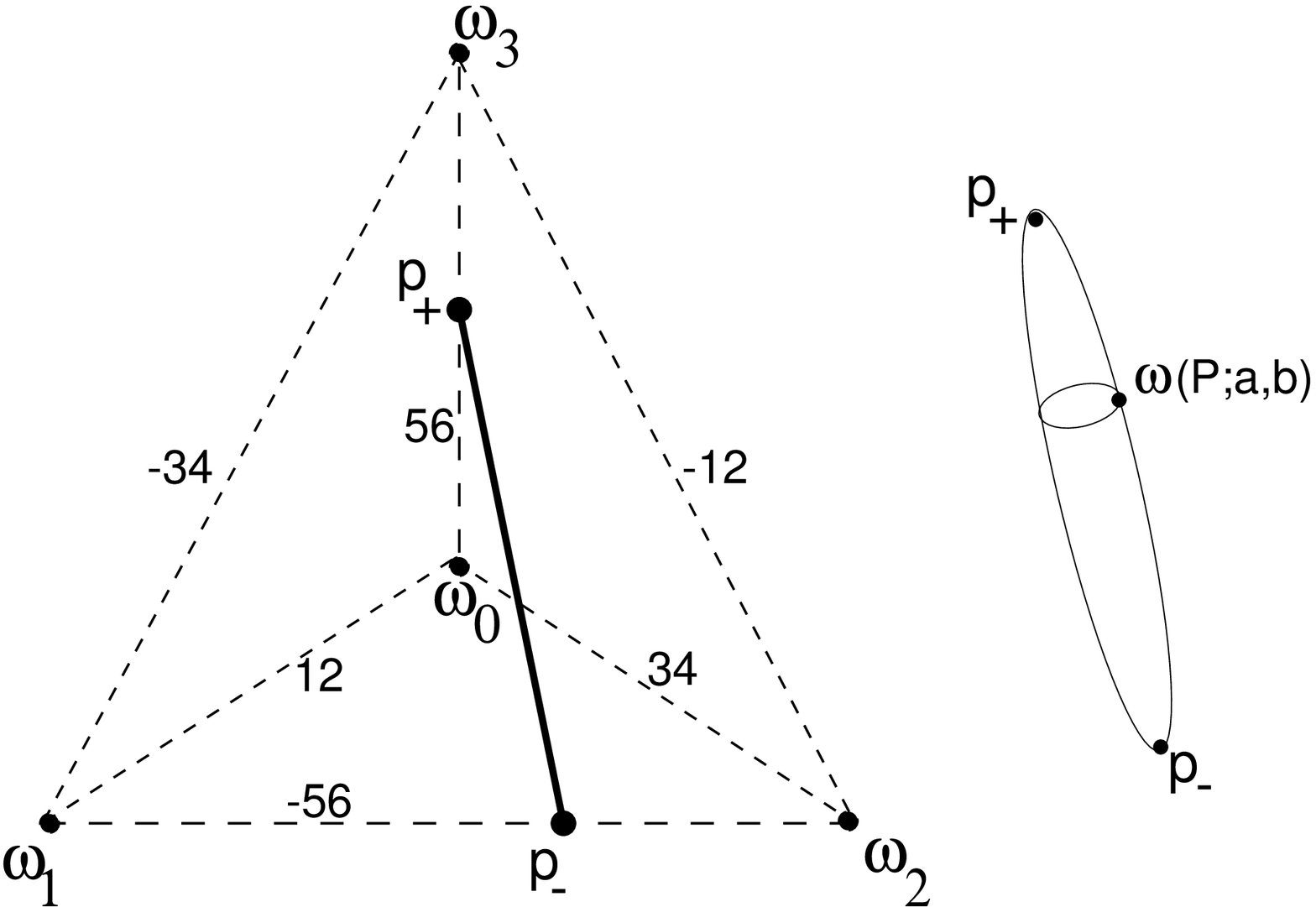}
\vspace{15pt}\vspace{30pt}

\nit{Example} Consider the subset \be\cS=\{\w(\m;0,1):\hbox{$\m_+$ has
$r=1$,\, $\m_-$ is arbitrary}\}\ee{S} of \rf{perp}.The condition $r=1$ implies
that $\m_+$ leaves fixed $e^{12}+e^{34}$, and so $\m(\w_0)=\w_0$. A point of
$\cS$ therefore lies on a `wire' joining $\w_0$ to an arbitrary point of
$\cE_{12}$, and is contained in the affine plane $\cF_3\setminus\cE_{12}
\cong\C^2$. Fix $e^5\we f^1+f^2\we e^6+f^{34}\in\cS$, and let $J$ be the
associated almost complex structure. The fact that $J$ lies in the face
opposite $\w_3$ implies that $J$ commutes with $J_3$ \cite{Ap}, so
\[f^2=Je^6=-JJ_3e^5=-J_3Je^5=J_3f^1\] is determined by $f^1$, which is a unit
vector in $\D$. It follows that $\cS$ is a 3-sphere, and the mapping
$\w(\m;0,1)\mapsto p_-$ is a Hopf fibration $\cS\to\cE_{12}\cong S^2$. It was
shown in \cite{AGS} that $\cS$ parametrizes the set of almost K\"ahler
structures on the Iwasawa manifold.

\subsection*{3. Detecting classes by exterior differentiation}\vs

Let $\nabla$ denote the Levi-Civita connection. The covariant derivative
$\nabla\w$ of the fundamental 2--form of an almost Hermitian manifold $M$ has
various symmetry properties. If $W$ is a real vector space of dimension
$2n\ge6$ with an almost complex structure $J$ and a real positive definite
inner product compatible with $J$, let $\W$ denote the subspace of
$W^*\ot\E^2W^*$ of the tensors with the same symmetry properties of
$\nabla\w$.

A decomposition \be\W=\W_1\op\W_2\op\W_3\op\W_4\ee{dec} is given in \cite{GH},
in which the four summands are invariant under the action of the unitary group
$U(n)$ and irreducible. The complete classification of almost Hermitian
manifolds into $2^4=16$ classes is obtained by taking $W$ to be the tangent
space. Ignoring the generic case, in which no component of $\nabla\w$
vanishes, we consider the subsets \be\ba{rcl}
\cZ_{ijk}&=&\{J\in\cZ:\nabla\w\in\W_i\op\W_j\op\W_k\},\y
\cZ_{ij}&=&\{J\in\cZ:\nabla\w\in\W_i\op\W_j\},\y
\cZ_i&=&\{J\in\cZ:\nabla\w\in\W_i \},\ea\ee{ijk} where $i,j,k$ are distinct
elements of $\{1,2,3,4\}$. Observe that, in this notation,
$\cZ_i\subseteq\cZ_{ij}\subseteq\cZ_{ijk}$, whereas $\cZ_i\cap\cZ_j$ is the
set of K\"ahler structures if $i\ne j$.

In particular,\par\sb\hspace{50pt}\bt{lcl} $\cZ_1$ &is the set of&
nearly-K\"ahler structures,\\ $\cZ_2$ && almost K\"ahler structures\\ $\cZ_3$
&& cosymplectic Hermitian structures,\\ $\cZ_4$ && locally conformal K\"ahler
structures.\et\par\sb\n The last statement makes use of the hypothesis $n\ge3$
in order that $d\w=\w\we\theta$, with $\theta$ a \textit{closed} 1--form.
However, we shall be more interested in determining the `maximal' classes
$\cZ_{234}$, $\cZ_{134}$, $\cZ_{124}$, $\cZ_{123}$, since all the others can
be obtained as intersections of these four. The only one of these that has a
special name is the class $\cZ_{123}$ of `semi-K\"ahler' or `cosymplectic'
structures, characterized by the condition $d\w\we\w^{n-2}=0$.

In \cite {FFS} the decomposition \rf{dec} is described in terms of
differential forms, and we review this approach. Given $J$, the complexified
cotangent space at any point $m$ of $M$ is given by \be T^*_m
M\otimes_\R\C=\E^{1,0}\op\E^{0,1}\ee{cot} where the two summands are the
$+1,\,-1$ eigenspaces of $J$ , respectively. This leads to the usual
decomposition of forms into types, whereby the $(p\+q)$th exterior power of
\rf{cot} contains a subspace $\E^{p,q}$ isomorphic to
$\E^p\E^{1,0}\otimes\E^q\E^{0,1}$. Both $\E^{p,q}\op\E^{q,p}, p\ne q$, and
$\E^{p,p}$ are complexifications of real vector spaces which are denoted by
$\ll\E^{p,q}\rr$ and $[\E^{p,p}]$ respectively in \cite{FFS}. Then

\nbf{Proposition 3} There are isomorphisms \[\W_1\cong\ll\E^{3,0}\rr,
\q\W_2\cong\ll V(2,1,0,\ldots,0)\rr,\q\W_3\cong\ll\E^{2,1}_0\rr,\q
\W_4\cong \ll\E^{1,0}\rr\cong T^*_mM.\] 

\n Here, $V(\la)$ denotes that irreducible complex representation of $U(n)$
with dominant weight $\la$. Also, $\E^{p,q}_0$ is the space of `effective'
forms, equivalently the Hermitian complement of the image of $\E^{p-1,q-1}$
under wedging with $\w$.\sb

It is well known that the covariant derivative $\nabla\w$ of the fundamental
2--form determines the Nijenhuis tensor of the almost complex structure $J$,
which is essentially the real part of the tensorial operator \be
d:\E^{1,0}\to\E^{0,2},\ee{Nij} and lies in $\W_1\op\W_2$. As observed in
\cite{FFS}, the skew-symmetric part of $\nabla\w$, which lies in
\[\E^3T^*_mM=\ll\E^{3,0}\rr\op\ll\E^{2,1}_0\rr\op
\ll\E^{1,0}\rr\cong\W_1\op\W_3\op\W_4,\] is proportional to $d\w$. These
algebraic properties lead to an easier way of computing the 16 Gary-Hervella
classes, which is especially fruitful in the case of $n=3$.

Suppose from now on that $M=G$ is a Lie group of real dimension 6. Let
$(\alpha,\beta,\gamma)$ be a basis of invariant $(1,0)$--forms, so that the
fundamental 2--form may be written \be\w=\frt12i(\alpha\we\bar
\alpha+\beta\we\bar\beta+\gamma\we\bar\gamma).\ee{aa} We omit the wedging
symbol `$\we$' for the remainder of this section.

\nbf{Lemma 1} The class $\cZ_{234}$ is characterized by the equation
\be(d\alpha)\beta\gamma\alpha{\bar\alpha}+(d\beta)\gamma\alpha
\beta{\bar\beta} +(d\gamma)\alpha\beta\gamma{\bar\gamma}=0.\ee{nK}

\nit{Proof} Recall that $\W_1\cong\ll\E^{3,0}\rr$. Since $\alpha\beta\gamma$
spans the space of invariant $(3,0)$--forms and $\nabla\w$ is \textit{real},
the component of $\nabla\w$ in $\W_1$ vanishes if and only if
$(d\w)\alpha\beta\gamma=0$. The result follows from \rf{aa}.\qed

\nbf{Lemma 2} The class $\cZ_{134}$ is characterized by the equations
\[\ba{c}
(d\alpha)\gamma\alpha\beta{\bar \beta}=0\\
(d\alpha)\alpha\beta\gamma{\bar \gamma}=0\\ 
(d\beta)\alpha\beta\gamma{\bar\gamma}=0\\
(d\beta)\beta\gamma\alpha{\bar \alpha}=0\\ 
(d\gamma)\beta\gamma\alpha{\bar \alpha}=0\\
(d\gamma)\gamma\alpha\beta{\bar\beta}=0\\[2pt]
(d\alpha)\beta\gamma\alpha{\bar\alpha}=(d\beta)\gamma\alpha\beta{\bar \beta}
=(d\gamma)\alpha\beta\gamma{\bar\gamma}.\ea\]\sb

\nit{Proof} The Nijenhuis tensor $N$ can be identified with the component of
$\nabla\w$ in $\W_1\op\W_2$, so we seek that component of $N$ belonging to
$\W_2$. With reference to \rf{Nij}, we may write \[N=((d\alpha)^{0,2},
(d\beta)^{0,2},(d\gamma)^{0,2}).\] Now, $(d\alpha)^{0,2}=0$ if and only if
\[(d\alpha)\beta\gamma\alpha\bar\alpha=(d\alpha)\gamma\alpha\beta\bar\beta
=(d\alpha)\alpha\beta\gamma\bar\gamma=0.\] The vanishing of $N$ would give 9
such complex--valued equations altogether, though a linear combination of
these is represented by Lemma~1. The proof is completed by the observation
that the difference of any two of the three terms in the last line of
equations is orthogonal to the left--hand side of \rf{nK}.\qed

\nbf{Lemma 3} The class $\cZ_{124}$ is characterized by the equations
$(d\w)\eta_i=0$ for $1\le i\le6$, where \[\ba{rcl} 
\eta_1&=&\alpha\beta{\bar\gamma}\\\eta_2&=&\beta\gamma{\bar\alpha}\\ 
\eta_3&=&\alpha\gamma{\bar \beta}\\
\eta_4&=&\alpha\beta{\bar\beta}-\alpha\gamma{\bar \gamma}\\
\eta_5&=&\beta\alpha{\bar\alpha}-\beta\gamma{\bar\gamma}\\
\eta_6&=&\gamma\alpha{\bar\alpha}-\gamma\beta{\bar\beta}.\ea\]\sb

\nit{Proof} Recall that $\W_3\cong \ll\E^{2,1}_0\rr$, and $(d\w)\eta_i$
represents the component of $d\w$ parallel to the effective $(2,1)$--form
$\bar\eta_i$.\qed

\n These equations can also be re--written in terms of $d\alpha$, $d\beta$,
$d\gamma$, but to no apparent advantage.

\nbf{Lemma 4} If the Lie group $G$ is \textit{nilpotent}, the class
$\cZ_{123}$ is characterized by the equations \vspace{-5pt}\be\ba{rcl}
(d\alpha)\w\w&=&0,\\ (d\beta)\w\w&=&0,\\(d\gamma) \w\w&=&0.\ea\ee{ww}

\nit{Proof} Given that $n=3$, $\cZ_{123}$ is determined by the condition \[
0=d(\w^2)=2\w\,d\w.\] The nilpotency condition implies that the exterior
derivative of any invariant 5--form is zero. This can be seen either by
considering the action of $d$ on a basis of the space of invariant 5-forms, or
by appealing to Nomizu's theorem \cite{N} on an associated nilmanifold $\Ga\bs
G$. Since \[\R\cong H^6(\Ga\bs G,\R)\cong\frac{\E^6\g^*}{d(\E^5\g^*)},\] it
cannot be the case that $d$ of any 5-form is non--zero. The identity
\[0=d(\w\w\alpha) =2(d\w)\w\alpha+\w\w (d\alpha)\] now gives the result. \qed

The equations \rf{ww} can be interpreted as a linear (rather than quadratic)
constraint on $\w\in\cZ$. The theory of the $*$ operator on an oriented inner
product space permits us to write \[\si\we(*\w)=g(\si,\w)\upsilon,\] where $g$
is the induced inner product on forms of the same degree, and $\upsilon$ is
the volume form \rf{up}. But (up to a universal constant) $*\w=\w\w$, and
\rf{ww} becomes \[0=g(\w,d\alpha)=g(\w,d\beta)=g(\w,d\gamma).\] It follows
that $M$ is cosymplectic if and only if $\w$ is orthogonal to the image of
$d$. In the notation \rf{perp}, if $(\si_1,\ldots,\si_k)$, $k=6-b_1$, is a
real basis of the image of \rf{d}, then
\be\cZ_{123}=\bigcap_{i=1}^k\,\ZZ{\si_i}.\ee{123}

\nbf{Proposition 4} Let $G$ be a 6--dimensional nilpotent Lie group. Then,
relative to any left--invariant metric, $\cZ_{123}\ne\emptyset$.

\nit{Proof} The nilpotency of $G$ means that there is a basis $(e^i)$ of
$\g^*$ such that $de^1=de^2=0$ and $de^k\in\E^2\al e^1,\ldots,e^{k-1}\ar$ for
$k\ge3$. Applying Gram--Schmidt, we may assume that this basis is
orthonormal. Since $e^4$ appears only in $de^5,\,de^6$, and $e^5$ only in
$de^6$, there exists an orthonormal basis $(f^1,f^2,f^3)$ of $\al
e^1,e^2,e^3\ar$ such that $e^4\we f^1$ and $e^5\we f^2$ are orthogonal to
$d(\g^*)$. Then $e^4\we f^1+e^5\we f^2\pm e^6\we f^3\in\cZ_{123}$.\qed

\subsection*{4.~Sixteen classes of the Iwasawa manifold}

In this section we shall compute the almost Hermitian structures on the
Iwasawa manifold by combining Proposition~2 of \S2 and Lemmas~1--4 of \S3. We
work in terms of the matrix $\m$ given by \rf{SO4},\rf{m+-}. In order to state
our main result, we regard a point of $\cZ\cong\CP^3$ as a fundamental 2-form
of the corresponding almost complex structure, and use freely the notation
$\cE_{ij}$, $\cF_i$ (introduced after \rf{vert}) and $\cS$ (from \rf{S}).

\nbf{Theorem 1} Let $M$ be the Iwasawa manifold endowed with the
metric \rf{g}. The classes of almost Hermitian structures defined in \rf{ijk}
are given by the following subsets of $\cZ$.\vspace{-2pt}\[\ba{rcl}
\cZ_{234}&=&\cF_3,\\ 
\cZ_{134}&=&\{\w_0\}\cup\cE_{12},\\
\cZ_{124}&=&\{\w_3\}\cup\cS,\\ 
\cZ_{123}&=&\cF_0\cup\cF_3.\ea\]\vspace{3pt}

\nbf{Corollary 1} The classes \rf{ijk} for $M$ are given by the lattice
\[\ba{lccc} \cZ_{123}\=&& \cF_0\cup\cF_3 &\\[24pt] \cZ_{23}\=\cZ_{234}\=&&&
\cF_3\\[12pt] \cZ_{12}\=\cZ_{124}\=& \{\w_3\}\cup\cS &&\\[8pt] &
&&\\[12pt] \cZ_2\=\cZ_{24}\=& \cS &&\\[8pt]
\cZ_3\=\cZ_{13}\=\cZ_{34}\=\cZ_{134}\=&&& \{\w_0\}\cup\cE_{12}\\[24pt]
\cZ_1\=\cZ_4\=\cZ_{14}\=&& \emptyset &\ea\]\vs

\vspace{-145pt}
\phantom.\hspace{200pt}\epsfxsize120pt\epsfbox[0 0 205 160]{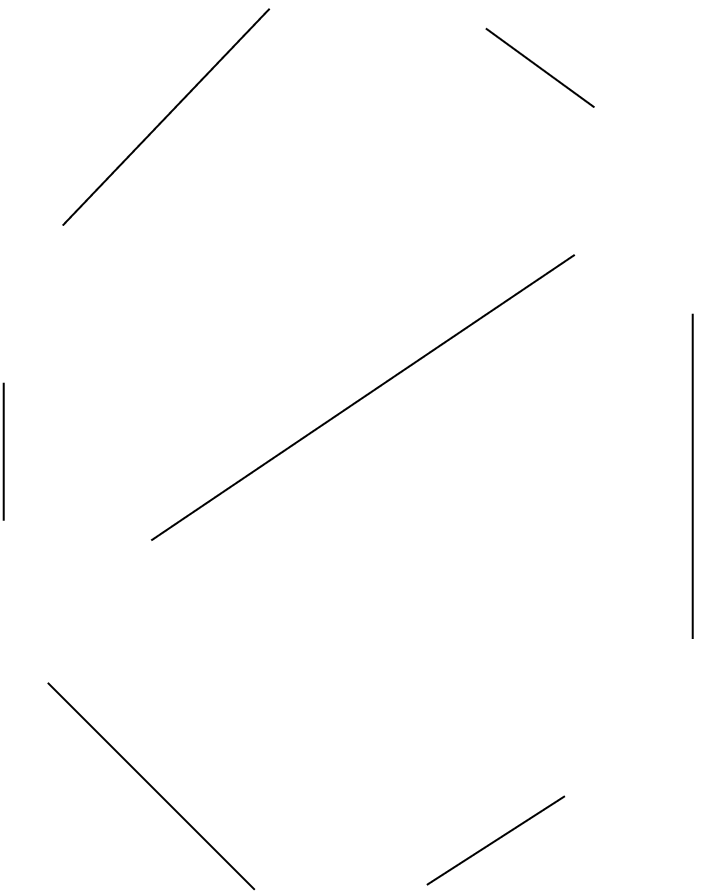}
\vspace{45pt}

\n The fact that $\cZ_{34}=\cZ_3$ means that all the Hermitian structures on
$M$ are automatically cosymplectic. We shall have more to say in \S6 regarding
other equalities between the classes.\sb

The elements \be\left\{\ba{rcl} \alpha&=&ae^6+be^1-ie^5,\y\beta&=&ae^1-be^6+
ie^2,\y \gamma&=&e^3+ie^4\ea\right.\ee{basis} are $(1,0)$--forms relative to
the fundamental 2--form \rf{w} in which $\m=\id$, and $\w(\id;a,b)$ is then
given by \rf{aa}. The problem is to determine $\m\in SO(4)\subset SO(6)$ and
$(a,b)\in S^1$ such that $\w(\m;a,b)=\m(\w(\id;a,b))$ satisfies the
appropriate condition. Now $\m\alpha,\m\beta,\m\gamma$ are $(1,0)$--forms
relative to $\w(\m;a,b)$, so (as an example) the first equation of Lemma~2
becomes \[(d(\m\alpha))(\m\gamma)(\m\alpha)(\m\beta)(\m\bar\beta)=0,\] or
equivalently \be(\dm\alpha)\gamma\alpha\beta\bar\beta=0,\ee{dm} where
$\dm=\m^{-1}\circ d \circ \m$. Our method then consists of applying
Lemmas~1,2,3,4 to \rf{basis}, but using the operator $\dm$ in place of $d$.

The fact that $\D=\ker d$ ensures that $d\circ \m=d$, and matters are even
simpler in the case of the complex Heisenberg group, for then
$\dm=\m^{-1}_+\circ d$. Indeed, \be d(\m\alpha)=a\,de^6-ide^5,\q
d(\m\beta)=-b\,de^6,\q d(\m\gamma)=0,\ee{dma} and we may assume that
\be\ba{rcl} \dm e^5 &=& s(e^{12}+e^{34}) +
v(e^{13}+e^{42})+y(e^{14}+e^{23}),\y \dm e^6 &=& t(e^{12}+e^{34}) +
w(e^{13}+e^{42})+z(e^{14}+e^{23}),\ea\ee{56} in accordance with \rf{rv} and
\rf{m+-} (noting that $\m_+^{-1}$ is represented by the transpose matrix).

With these preliminaries, we may proceed to the proof of Theorem~1.

\nbf{The class $\cZ_{234}$} From \rf{basis}, \[\ba{l}\alpha\beta
\gamma=(e^{136}-ae^{246}-be^{124}+ae^{145}+e^{235}-be^{456})
\\\ph{mmmmmm}+i(e^{146}+ae^{236}+be^{123}-ae^{135}+be^{356}+e^{245}).\ea\]
The equation defining $\cZ_{234}$ can now be obtained from Lemma~1, with the
aid of \rf{56}. In the case of the Iwasawa manifold, it is \[ (a+1)[(y + w)
+i(z-v)]=0.\]

\n The solutions are

\n(1i) $a=-1$, giving the entire edge $\cE_{12}$ of the tetrahedron.

\n(1ii) $\{y=-w,\,z=v\}$, giving \be\m_+=\left(\ba{ccc}1&0&0\\0&v&-w\\0&w&v\ea
\!\right),\q v^2+w^2=1,\ee{1ii} bearing in mind that the $3\times3$ matrix
must be orthogonal with determinant 1. For arbitrary $\m\in SO(4)$,
$\w(\m;1,0)$ generates the edge $\cE_{03}$, but here, $e^{12}+e^{34}$ is fixed
by $\m$, so $\w(\m;1,0)=\w_0$. Since $\w(\m;-1,0)$ is still free to move on
$\cE_{12}$, we obtain the whole face $\cF_3$ (compare \rf{S}) as $a,b$ vary.

\n Note that (1i) is now subsumed under (1ii).

\nbf{The class $\cZ_{134}$} Combining \rf{dma} with Lemma~2 gives the
respective equations \[\left\{\ba{l} b[(v -a z) +i(a w + y)]=0,\y
b[at-is]=0,\y b^2t=0,\y b(1+ a)[z-iw]=0,\y (a + 1)[(-2 a z + v + z) + i(2 a w
- w + y)]=0,\y (a + 1)[(v - a z) +i(a w +y)]=0.\ea\right.\]

\n Possible solutions are:

\n(2i) $\{b=0,\,a=-1\}$, giving the edge $\cE_{12}$. 

\n(2ii) $\{b=0,\,a=1,\,v=z,\,w=-y\}$, giving the vertex $\w_0$ as a special
case of \rf{1ii}.

\n(2iii) If $b\ne0$ then the middle two equations give $t=w=z=0$. This would
imply that $\m_+$ is singular, which is impossible.

\nbf{The class $\cZ_{124}$} For this, one first computes the $\eta_i$ in
Lemma~3. For
example,\[\ba{l}\eta_1=e^{136}+ae^{246}+be^{124}-ae^{145}+be^{456}+e^{235}
\\\ph{mmmmmmmm}+i(-e^{146}+ae^{236}+be^{123}-ae^{135}+be^{356}-e^{245}).
\ea\]

\n Lemma~3 now gives a total of 11 real equations: \[\left\{\ba{l} (a-1)[(y +
w)+i(v-z)]=0,\y (a + 1)[(2 a w - y - w)+i(2 a z - z + v)]=0,\y (a - 1)[(2 a w +
y + w)+i(2 a z + z - v)]=0,\y (a-1)[s +iat]=0,\y ib(a-1)t=0,\y ab(
w+iz)=0.\ea\right.\]
The solutions are

\n(3i) $\{a=1,\,b=0,\,y=w,\,v=-z\}$, giving \be \m_+=\left(\ba{ccc}-1&0&0\\
0&v&w\\0&w&-v\ea\right ),\ee{3i} and the vertex $\w_3$.

\n(3ii) $b\ne0$, forcing $s=t=0$ and $y=-w,\,v=z$. If $a\ne0$ we obtain
$w=z=0$, which is impossible on top of $t=0$. Whence, $a=0$, $b=\pm1$, and
$\m_+$ satisfies \rf{1ii}. It follows that $\w(\m;0,\pm1)$ describes the
3-sphere \rf{S}.

\nbf{The class $\cZ_{123}$} We know from Corollary~1 that this is the
intersection of $\ZZ{e^{13}\!+\!e^{42}}$ and $\ZZ{e^{14}\!+\!e^{23}}$, though
we shall use the $SO(4)$ action to determine it geometrically. Observe that
\[\ba{rcl} d\w\we\w &=& de^5\we(e^{126}+ae^{346}+be^{134})-de^6\we(e^{125}+
ae^{345}+be^{234}) \y &=& (a + 1)s\,e^{12346}-(a + 1)t\,e^{12345}.\ea\] This
can vanish in one of two ways:

\n(4i) $a=-1$, giving the edge $\cE_{12}$.

\n(4ii) $s=t=0$, which reduces $\m_+$ to one of \rf{1ii},\,\rf{3i}. In the
former case ($r=1$), the vertex $\w_0$ is fixed on $\cE_{03}$ and we obtain
the face $\cF_3$, as explained in (1ii). In the latter case ($r=-1$),
$\m(\w_0)=\w_3$, and we obtain the face $\cF_0$ opposite to $\w_0$.

\n In conclusion, $\cZ_{123} =\cF_0\cup\cF_3$.

\sb The proof of Theorem~1 is now complete.

\subsection*{5.~An example with discrete Hermitian structures}

We have explained (after \rf{dm}) the importance of the condition that the
4--dimensional $SO(4)$ module $\D$ used to define $\w(\m;a,b)$ lie in the
kernel of $d\colon\g^*\to\E^2\g^*$. According to the classification, there are
7 isomorphism classes of nilpotent Lie algebras that, in common with $\gH$,
have first Betti number $b_1$ (the dimension of the kernel of \rf{d}) equal to
4. Of these, only three are irreducible, namely $\g_1=\gH$ given by \rf{rv},
another 2--step Lie algebra $\g_2$, and a 3--step algebra $\g_3$ (see for
example \cite{S}). There is no particular reason to restrict to irreducible
algebras, except that one would expect the reducible cases to be easier to
describe.

The structure equations \be\left\{\ba{l}de^i=0,\q 1\le i\le4\y de^5=
e^{12},\y de^6=e^{14}+e^{23}.\ea\right.\ee{g2} of $\g_2$ are very similar to
\rf{rv}. In applying the methods of \S4, the real only difference is that now
\be\ba{rcl} \dm e^5&=&
\m^{-1}\Big(\frac12(e^{12}+e^{34})+\frac12(e^{12}-e^{34})\Big)\\[6pt]
&=&\frac12\Big[r(e^{12}+e^{34})+u(e^{13}+e^{42})+x(e^{14}+e^{23})\y
&&\ph{mmmmmm}+r'(e^{12}-e^{34})+u'(e^{13}-e^{42})+x'(e^{14}-e^{23})\Big]
\ea\ee{12} (compare \rf{56}).

As usual, we consider the $e^i$ as invariant forms on an associated Lie group
$G_2$, or nilmanifold $M_2$, and we choose the metric $g$ (as in \rf{g}) that
renders them orthonormal.

\nbf{Theorem 2} Let $M$ be a nilmanifold associated to \rf{g2}, with the
metric $g$. The classes of invariant almost Hermitian structures defined in
\rf{ijk} satisfy \vspace{-2pt}\[\ba{rcl} \cZ_{234}&\supset&\{\w_1,\w_2,
\ep(\frac{2\pi}3),\ep(-\frac{2\pi}3)\}\cup\CS,\y\cZ_{134}&=&\{\w_1,\w_2,
\ep(\frac{2\pi}3),\ep(-\frac{2\pi}3)\},\y\cZ_{124}&=& \CS,
\y\cZ_{123}&\supset&\{\ep(0),\ep(\pi)\}\cup\CS\cup\CS'\cup\CS_{02}\cup
\CS_{13}\cup \CS_{12},\ea\] (notation as in \rf{equat},\rf{eth}), where
$\CS,\CS'$ are other circles and the unions are disjoint.\mb

Provided we use \rf{12}, the computations required for the proof of Theorem 2
are very similar to those of \S4, and we omit verification of facts (2i) and
(3i) below that are used to determine the `larger classes'.

\nbf{The class $\cZ_{234}$} Lemma~1 gives \[
(a+1)(x+2w)+(a-1)x'-i[(a+1)(u-2z)+(a-1)u']=0.\]

\n(1i) When $a=-1$ we get $x'=0=u'$ giving the two points $\w_1,\w_2$.

\n(1ii) When $a=1$ we get $x=-2w$ and $u=2z$. This implies that one of $r,t$
vanishes, and the only possibility is that $r=0$ and $4z^2+4w^2=1$. This
forces $t=\pm\frac12\sqrt3$ and $s=-\frac12$ (the sign of $s$ is fixed by the
condition that $\det\m_+=1$), giving solutions
\[\ba{rcl}\ep(\frac{2\pi}3)&=&-\frac12(e^{13}\+e^{42})+\frac{\sqrt3}2
(e^{14}\+e^{23})+e^{56},\y \ep(-\frac{2\pi}3)&=&-\frac12(e^{13}\+e^{42})
-\frac{\sqrt3}2(e^{14}\+e^{23})+e^{56}.\ea\]

\n(1iii) When $a=0$ we get $x'=x+2w$ and $u'=u-2z$, that include the solutions
in (3ii) below.

\nbf{The class $\cZ_{134}$} Lemma~2 gives 

\n(2i) If $b\ne0$ then $\w(\m;a,b)$ cannot satisfy $W_2=0$.

\n(2ii) If $a=1$ the equations are the same as in (1ii).

\n(2iii) If $a=-1$ we get $u'=0=x'$, giving $\w_1,\w_2$.

\nbf{The class $\cZ_{124}$} 

\n(3i) If $a\ne0$ then $\w(\m;a,b)$ cannot satisfy $W_3=0$.

\n(3ii) If $a=0$, Lemma~3 gives \be u'=u-2z,\q x'=x+2w,\q r'=r,\q t=0.\ee{u'u}
Let $\xi$ denote the unit vector $(r,u,x)$. Then we see that
\[|\xi-\xi'|^2=(r-r')^2+(u-u')^2+(x-x')^2=4(w^2+z^2)=4.\] This forces
$\xi'=-\xi$, whence $r=0$ and \[\ba{l}
\m^{-1}(e^{12})=ue^{42}+xe^{23}=e^2\we(xe^3-ue^4),\y
\m^{-1}(e^{34})=xe^{14}+ue^{13}=e^1\we(xe^4+ue^3).\ea\] Since the top row of
$\m_+$ is determined by $P^{-1}(e^{12}+e^{34})$, and $\det\m_+=1$, we obtain
\[\m_+=\left(\ba{ccc}0&u&x\\1&0&0\\0&x&-u\ea\right)\q\hbox{or}\q
\left(\ba{ccc}0&u&x\\-1&0&0\\0&-x&u\ea\right).\] However, the first matrix
gives $x'=3x$, $u'=3u$ and $|\xi|\ge3$, which is impossible. 

In computing $\w(\m;0,1)=P(\varpi_0)$, we are free to apply an element of
$SO(2)$ that rotates $e^3,e^4$ so that $x=1$ and $u=0$ (see \rf{varpi},
\rf{rot}). It follows that $\m^{-1}\al e^1,e^2\ar=\al e^2,e^3\ar$ and
$\m^{-1}\al e^3,e^4\ar=\al e^1,e^4\ar$, and \[\m^{-1}:\q\ba{ll} e^1\mapsto
Ae^2+Be^3\q&\q e^2\mapsto-Be^2+Ae^3,\y e^3\mapsto Ae^1-Be^4\q&\q e^4\mapsto\
Be^1+Ae^4,\ea\] where $A^2+B^2=1$, and one can check that \rf{u'u} holds. The
inverse mapping is given by \[\m:\q\ba{c} e^1\mapsto Ae^3+Be^4\qquad
e^2\mapsto Ae^1-Be^2,\y e^{34}\mapsto(Be^1+Ae^2)(-Be^3+Ae^4),\ea\] and it
follows that $\cZ_{124}$ is the circle \be\CS=\{-(Ae^3+Be^4)\we
e^5+(Ae^1-Be^2)\we e^6 +(Be^1+Ae^2)\we(-Be^3+Ae^4)\}\ee{CS} in the `middle' of
the tetrahedron.

\nbf{The class $\cZ_{123}$} In order to visualize this using Corollary~1, it
is helpful to realize that the hypersurface $\ZZ{e^{12}}$ is `suspended' half
way between the edges $\cE_{01}=\EE{e^{12}}$, $\cE_{23}=\EE{-e^{12}}$ of the
tetrahedron (labelled `$12$' and `$-12$' in the figure in \S2), with which it
has empty intersection. Modulo a change of basis, this description follows
from \rf{perp}.

The additional orthogonality to $e^{14}+e^{23}$ implies the following:

\n(4i) $\cZ_{123}\cap\cE_{03}$ consists of the points
\[\ep(0)=e^{13}+e^{42}+e^{56},\q \ep(\pi)=-e^{13}-e^{42}+e^{56}.\]

\n(4ii) $\cZ_{123}$ contains the full equators \rf{equat} in
$\cE_{12},\cE_{02}, \cE_{13}$.

\n(4iii) Any point on the circle \rf{CS} is obviously orthogonal to both
$e^{12}$ and $e^{14}+e^{23}$. The same is true for the circle
\[\CS'=\{(Ae^3+Be^4)\we e^5+(Ae^1-Be^2)\we e^6-(Be^1+Ae^2)
\we(-Be^3+Ae^4)\}.\]

This completes the proof of Theorem~2, which provides enough information to
determine most of the classes $\cZ_i$, $\cZ_{ij}$.

\nbf{Corollary~2} The set $\cZ_{34}$ of invariant Hermitian structures on
$(M_2,g)$ consists of four distinct points on $\CS_{03}$, and the set $\cZ_2$
of compatible symplectic structures is the circle $\CS$.\mb

According to \cite[Theorem~3.3]{S}, \textit{any} invariant complex structure
$J$ on $M_2$ has the property that $\D=\al e^1,e^2,e^3,e^4\ar$ is
$J$-invariant. Let $\E^{p,q}\D=\E^{p,q}\g^*\cap\D_\C$. Any $(1,0)$-form not in
$\E^{1,0}D$ will be a multiple of $\alpha=e^5+ce^6+\beta$, with $c\in\C$ and
$\beta\in\E^{1,0}\D$. The classification of complex structures then amounts to
determining which almost complex structures on $\D$ have the property that
\[d\alpha=e^{12}+c(e^{14}+e^{23})\in\E^{1,1}\D\op\E^{2,0}\D.\] This becomes a
more subtle problem with the imposition of an arbitrary metric.\sb

\nit{Example} Let $g'=\frac14(e^1\ot e^1+e^2\ot e^2)+\sum_{i=3}^6e^i\ot e^i$,
so that $\{\frac12e^1,\frac12e^2,e^3,e^4,e^5,e^6\}$ is orthonormal relative to
$g'$. We denote the space of $g'$-orthogonal almost complex structures by
$\cZ'$. Since any complex structure $J\in\cZ'$ satisfies $J(\D)=\D$, it must
be that $Je^5=\pm e^6$ and $J$ lies in $\cE'_{12}$ ($c=i$) or $\cE'_{03}$
($c=-i$). Two obvious candidates in $\cE_{12}'$, namely \[\ba{l}
\w_1'=\;\frt14e^{12}-e^{34}-e^{56},\y \w_2'=-\frt14e^{12}+e^{34}-e^{56},\ea\]
are the counterparts to $\w_1,\w_2\in\cZ_{34}$. But this time, the equation
\[d(e^5+ie^6)=e^{12}+i(e^{14}+e^{23})=e^1 (\frt12e^2+ie^4)+
e^2(-\frt12e^1+ie^3)\] leads to a \textit{unique} solution \[
\ep'(\pi)=-\frt12e^{13}-\frt12e^{42}+e^{56}\in\CS_{03}'.\] Thus, $\cZ_{34}'=
\{\w_1',\w_2',\ep'(\pi)\}$ consists of exactly \textit{three} points.

This example is significant, because it shows that the homotopy class of the
set of Hermitian structures is not independent of the choice of metric. One
would expect the set of invariant Hermitian structures on $M_2$ to be discrete
relative to any inner product on $\g_2$. In this connection, the maximum
number of isolated orthogonal complex structures that can co--exist on
6-manifold will depend on properties of its Weyl conformal curvature tensor.

\subsection*{6.~Treatment of the 3--step case}

A greater contrast with the calculations of \S4 is provided by $\g_3$, whose
structure equations are \be\left\{\ba{l}de^i=0,\q 1\le i\le4\y de^5= e^{12},\y
de^6=e^{15}+e^{34}.\ea\right.\ee{g3} This example is of special interest,
since the associated group $G_3$ or nilmanifold $M_3$ possesses \textit{no}
invariant complex \textit{or} symplectic forms, relative to \textit{any}
metric \cite{S}. We shall detect and generalize this non--existence for the
structures compatible with the standard metric \rf{g}.

Equation \rf{12} is again applicable. In order to compute $\dm e^6$, recall
that $\w(\m;a,b)=\w(\m Q;a,b)$ if $Q$ belongs to the subgroup $SO(2)$ of
$SO(4)$ which fixes $e^1,e^2$, and rotates $e^3,e^4$. By left--multiplying
$\m$ by a suitable $Q$, we may suppose that \[\m^{-1}(e^1)=\suml{i=1}3\la_i
e^i,\] where $\la_1^2+\la_2^2+\la_3^2=1$. With this choice of $\m$,
\[\ba{rcl}\dm e^6&=&\m^{-1}(e^1)\we e^5+\m^{-1} (e^{34})\y
&=&(\la_1e^{15}+\la_2e^{25}+\la_3e^{35})+\frac12\Big[r(e^{12}+
e^{34})+u(e^{13}+e^{42})+x(e^{14}+e^{23})\y &&\hspace{140pt}-r'(e^{12}-e^{34})
-u'(e^{13}-e^{42})-x'(e^{14}-e^{23})\Big].\ea\]

The presence of the $\la_i$ complicates the situation, though we shall give a
partial description of classes of almost Hermitian structures in terms of the
generalized edges \rf{EE}.     

\nbf{The class $\cZ_{234}$} According to Lemma~1, $\w(\m;a,b)$ has $W_1=0$ if
and only if \[\ba{l} x(1+a)-x'(1-a)+u(1+a)-u'(1-a)\y\ph{ooooo}+
i\Big[x(1+a)-x'(1-a)-u(1+a)+ u'(1-a) + 2b\la_3\Big]=0.\ea\] In the special
case $\la_3=0$, the equations become \be\left\{\ba{l} x(1+a)=x'(1-a),\y
u(1+a)=u'(1-a).\ea\right.\ee{xu} From \rf{12}, \[\m^{-1}(e^{12})=(\la_1e^1+
\la_2e^2)\we \m^{-1}(e^2)=\frt12(r+r')e^{12}+\frt12 (r-r')e^{34}+\cdots\]
Since the left-hand side has no term in $e^{34}$, we obtain $r=r'$. Thus,
$u^2+x^2=(u')^2+(x')^2$, and \rf{xu} gives two cases:

\n(1i) $u=x=u'=x'=0$.  Either $r=r'=1$, giving the edge $\cE_{02}$, or
$r=r'=-1$, giving the edge $\cE_{13}$.

\n(1ii) $a=0$. We have $r=r'$, $u=u'$ and $x=x'$. From \rf{12} we obtain \be
\m^{-1}(e^{12})=re^{12}+ue^{13}+xe^{14}=e^1\we(re^2+ue^3+xe^4).\ee{simple} It
follows that $\m^{-1}(e^1)$ is a linear combination of $e^1$ and
$re^2+ue^3+xe^4$, but from above this can only happen if $u=x=0$ (case (1i))
or $\m^{-1}(e^1)=\pm e^1$. In the latter case, the fundamental 2--form belongs
to the disjoint union $\EE{e^{15}}\,\cup\,\EE{-e^{15}}$ of two spheres,
parametrized by $(r,u,x)$ and the choice of sign.\sb

In view of these results, we perform a rotation of $90^{\mathrm{o}}$ in the
$\al e^2,e^5\ar$ plane to transform \rf{vert} into \be\ba{rcl} \bu\varpi_0
&=&-e^{15}+e^{26}+e^{34},\\\bu\varpi_1&=&-e^{15}-e^{26}-e^{34},
\\\bu\varpi_2&=&e^{15}-e^{26}+e^{34},\\\bu\varpi_3&=&e^{15}+e^{26}-e^{34},
\ea\ee{varvert} extending the notation \rf{varpi}. We now state

\nbf{Theorem 3} Let $M_3$ be a nilmanifold associated to \rf{g3}, with a
standard metric \rf{g}. The classes of almost Hermitian structures defined in
\rf{ijk} are given by\vspace{-2pt}\[\ba{rcl}
\cZ_{234}&\supset&\cE_{02}\cup\cE_{13}\cup\EE{e^{15}}\cup\EE{-e^{15}},\\
\cZ_{134}&=&\>\emptyset,\\\cZ_{124}&=&\{\varpi_1,\varpi_2\},\\
\cZ_{123}&\supset&\CS_{02}\cup \CS_{13}\cup\EE{e^{26}},\ea\] where 
$\varpi_1\in\CS_{13}$ and $\varpi_2\in\CS_{02}$.\mb

The first inclusion follows from (1i),(1ii), and is certainly strict. We
proceed to investigate the class $\cZ_{134}$ that in general contains the set
$\cC=\cZ_{34}$ of Hermitian structures.

\nbf{The class $\cZ_{134}$} Lemma~2 gives
\[\left\{\ba{l}
2a\la_3-(u+u')b+(x-x')ab-i[(x+x')b+(u-u')ab]=0,\y 
a\,[2a\la_1+b(r-r')]+i[2a\la_2-b(r+r')]=0,\y 
b\,[2a\la_1+b(r-r')]-2i\la_2b=0,\y 
b\,[(1+a)x+(1-a)x']-ib[(1+a)u+(1-a)u']=0,\y 
-(1+a)u+(1-a)u'+(a+a^2-b^2)x+(a-a^2+b^2)x'-2\la_3b\\
\ph{oo}-i\Big[(1+a)x-(1-a)x'+(a+a^2-b^2)u+(a-a^2+b^2)u'\Big]=0,\y 
-(1+a)u+(1-a)u'+(a+a^2)x+(a-a^2)x'\\
\ph{oo}-i\Big[(1+a)x-(1-a)x'+(a+a^2)u+(a-a^2)u'\Big]=0.\ea\right.\]

\n We obtain the following cases.

\n(2i) $a=0,\,b=\pm1$ implies $r=r'=u=u'=x=x'=0$, which is
impossible since $\m_\pm$ must be non--singular.

\n(2ii) $b=0,\,a=\pm1$ implies that $\la_i=0$ for all $i$, which is
equally absurd. 

\n(2iii) $ab\ne0$ implies that $\la_2=\la_3=0$ and $u=u'=x=x'=0$.
Another consequence of the equations is that $r=-r'$, but this contradicts
Lemma~5.\sb

\n Thus, $\cZ_{134}$ is indeed empty. 

\nbf{The class $\cZ_{124}$} Lemma~3 gives
\[\left\{\ba{l}
(a-1)(u+x)+(a+1)(u'+x')\\
\ph{oo}+i[(a-1)(u-x)+(a+1)(u'-x')-2\la_3b]=0,\y 
-(a+1)x-(a-1)x'+(a+2a^2-1)u+(a-2a^2+1)u'\\
\ph{oo}+i[(a+1)u+(a-1)u'+(a+2a^2-1)x+(a-2a^2+1)x'+2\la_3b]=0,\y 
(a-1)x+(a+1)x'-(a-2a^2+1)u-(a+2a^2-1)u'\\
\ph{oo}-i[(a-1)u+(a+1)u'+(a-2a^2+1)x+(a+2a^2-1)x'+2\la_3b]=0,\y 
(a-1)r+(a+1)r'-ia[2\la_1b-(a-1)r+(a+1)r']=0\y 
b\,[2\la_1b-r(a-1)+r'(a+1)]=0\y 
ab\,[(u-u')+i(x-x')]=0.\ea\right.\]

\n Possible solutions are

\n(3i) $a=0$, $b=\pm1$ implies that $r=r'=\mp\la_1$, $u=u'$, $x=x'$, and
$\la_3=0$. From \rf{simple}, we deduce that $r=\pm1$ \textit{and} $u=x=0$. All
solutions lie on $\cE_{02}$ (if $r=r'=1$) or $\cE_{13}$ (if
$r=r'=-1$). However, the condition that $\m(e^1)=-bre^1$ reduces solutions to
$\varpi_1,\varpi_2$.

\n(3ii) $b=0$, $a=\pm 1$ implies that $r=u=x=r'=u'=x'=0$, which is impossible.

\n(3iii) $ab\ne0$ implies $u=x=u'=x'=0$, whence $r^2=(r')^2=1$. But this
contradicts $(1-a)r=(1+a)r'$, which is another consequence of the equations.

\nbf{The class $\cZ_{123}$} According to \rf{123}, this is a singular
intersection of $\ZZ{e^{12}}$ and $\ZZ{e^{15}\!+\!e^{34}}$. Lemma~4 gives
\[\left\{\ba{l} (a+1)r+(a-1)r'=0\\(a+1)r+r'(1-a)=2b\la_1,\ea\right.\] and thus
$\la_1b=r(a+1)=r'(1-a)$. We can use this to describe some subsets of
$\cZ_{123}$:

\n(4i) $b=0$ gives $r=0$ and solution set $\CS_{03}$, or $r'=0$ and $\CS_{12}$.

\n(4ii) $a=0$ implies $r=r'=\la_1$, and this gives $\EE{e^{26}}$.\mb

The subset $\cZ_{123}$ intersects the edges $\bigcup\cE_{ij}$ of the original
tetrahedron in the circles $\CS_{02},\CS_{13}$ and the points
$\varpi_0,\varpi_3\in\EE{e^{26}}$. It intersects each face $\cF_i$ in a `cone'
joining one of these two points to one of the two circles. On the other hand,
Theorem 3 is saying is that $\cZ_{234}$ contains four of the six edges
\[\cE_{02}=\EE{e^{34}},\q\cE_{13}=\EE{-e^{34}},\q\EE{\pm e^{15}},\q \EE{\pm
e^{26}}\] of the new tetrahedron with vertices \rf{varvert}. The new edge
$\EE{e^{26}}$ is \textit{not} one of these four, and is disjoint from
$\varpi_1,\varpi_2$, that one can check do not lie in $\cZ_{123}$. It follows
that the class $\cZ_2$ of symplectic structures is empty, and moreover

\nbf{Corollary 3} The following classes \rf{ijk} for $\!(M_3,g)$ are empty:
$\!\cZ_1,\cZ_2,\cZ_3,\cZ_4$, $\cZ_{12},\cZ_{13}$,
$\cZ_{14},\cZ_{34},\cZ_{134}$.\mb

As already remarked, the vanishing of $\cC=\cZ_{34}$ and $\cS=\cZ_2$ does not
depend on the choice of metric on $M_3$; it would be interesting to know
whether this is true of the other classes listed above. A more positive
feature of $M_3$ is that (at least with respect to $g$) $\cZ_{24}\setminus
\cZ_2$ is non--empty. Indeed, the equations \[ d\varpi_1=-\varpi_1\we e_2,\q
d\varpi_2=\varpi_2\we e_2\] show that $\cZ_{24}$ contains $\varpi_1,\varpi_2$,
and Theorem~3 implies that it contains no other points.\sb

The examples studied in \S4,\,\S5,\,\S6 all have $\cZ_{14}=\emptyset$ (see for
example Corollary~1), and it is natural to ask to what extent this is a
general phenomenon. In particular, we conjecture that the classes
$\cZ_1,\cZ_4$ are both empty for any Riemannian metric on any compact
6--dimensional nilmanifold other than a torus. This would generalize the
non--existence of a K\"ahler metric on such manifolds.

\renewcommand{\thebibliography}{\list{{\bf\arabic{enumi}.\hfil}}
{\settowidth\labelwidth{18pt}\leftmargin\labelwidth\advance
\leftmargin\labelsep\usecounter{enumi}}\def\newblock{\hskip.05em} \sloppy
\sfcode`\.=1000\relax} \newcommand{\bi}{\vspace{-4pt}\bibitem}

\enddocument